\documentclass[11pt,a4paper,final,dvips]{article}
\usepackage[cp1251]{inputenc}
\usepackage[russian,english]{babel}
\usepackage[T2A]{fontenc}

\usepackage{amsmath}
\usepackage{amsfonts}
\usepackage{amssymb}
\usepackage{graphicx}

\begin{document}

\begin{center}

  \textbf {On lines in a triangle tangent to a conic}

  Dolgirev Pavel

\end{center}

 \textbf {Abstract.} 
We present generalizations of theorems on {\it Kypert’s construction} and on 
{\it 2nd Morley’s Centre}. 
Most of our proofs are synthetic. 

\smallskip
\begin{center}
Introduction and main results
\end{center}

Let us introduce necessary notation.

Let us fix $\triangle ABC$. Draw lines $l_{A}$, $l_{B}$ and $l_{C}$ through vertices $A$, $B$ and $C$, respectively. Lines $l'_{A}$, $l'_{B}$ and $l'_{C}$ are symmetric to $l_{A}$, $l_{B}$ and $l_{C}$ with respect to bisecting lines of $\triangle ABC$. Let $X$, $Y$ and $Z$ be the intersection points of $l_{B}$ and $l'_{C}$, $l_{C}$ and $l'_{A}$, $l_{A}$ and $l'_{B}$, respectively. Let $X'$, $Y'$ and $Z'$ be the intersection points of $l'_{B}$ and $l_{C}$, $l'_{C}$ and $l_{A}$, $l'_{A}$ and $l_{B}$, respectively. Let $A_1$, $A'_1$, $B_1$, $B'_1$, $C_1$ and $C'_1$ be the intersection points of lines $l_{A}$, $l'_{A}$, $l_{B}$, $l'_{B}$ and $l_{C}$, $l'_{C}$ with sidelines of $\triangle ABC$.

\includegraphics{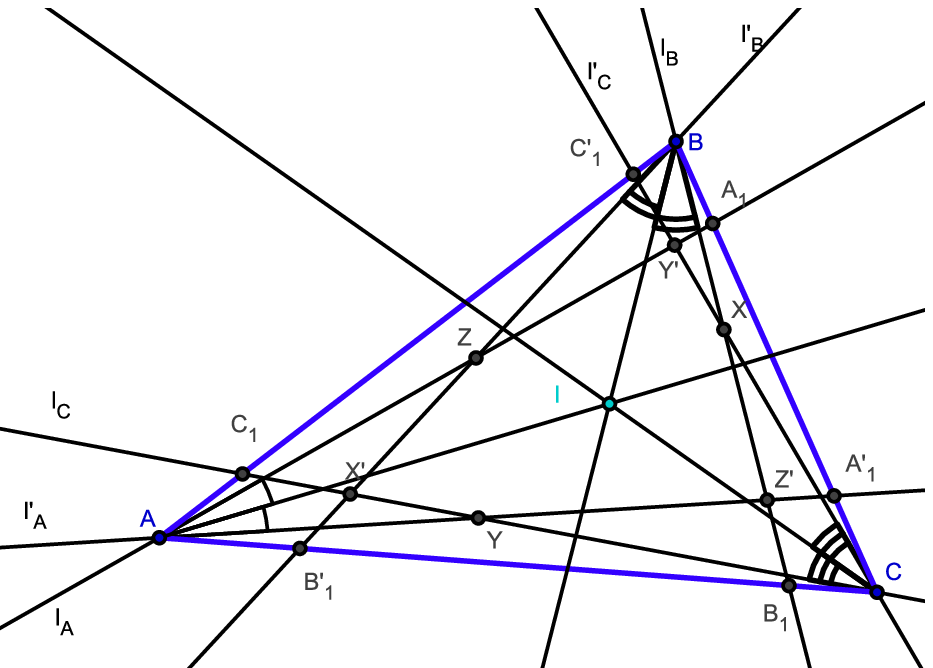}

\textbf{Theorem 1.} {\it Lines $AX$, $BY$ and $CZ$ intersect at a point.}

\includegraphics{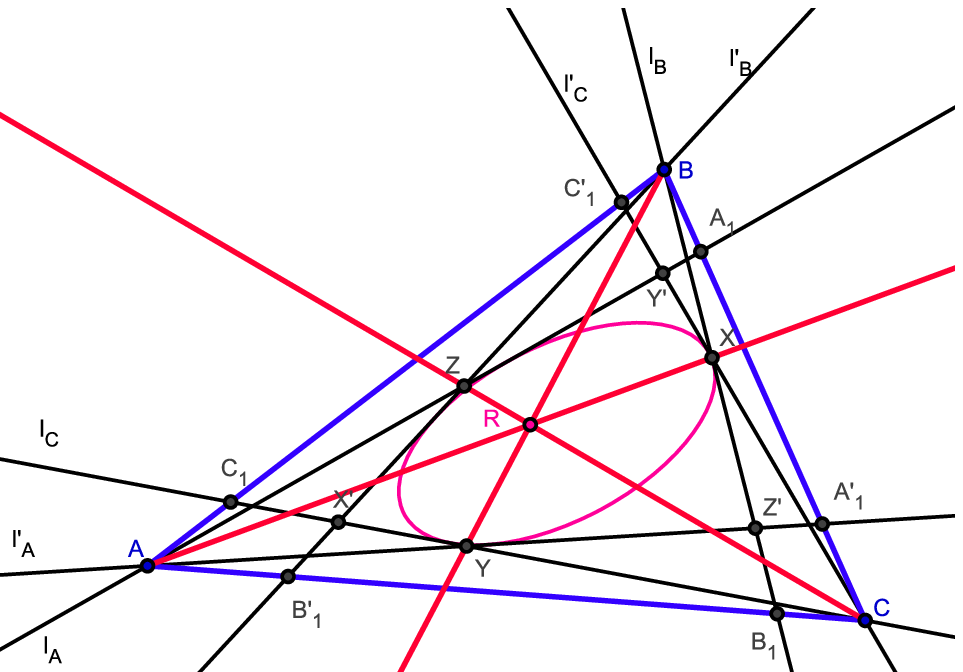}

\smallskip
\textbf{Corollary.} {\it Lines $l_{A}$, $l'_{A}$, $l_{B}$, $l'_{B}$, $l_{C}$ and $l'_{C}$ are tangent to a conic.}

\smallskip
\textbf{Theorem 2.} {\it Let $H_A$ be the intersection point of $XH_A$ and the line perpendicular to $BC$ and passing through $X$. Define $H_B$ and $H_C$ analogously.
Then lines $AH_A$, $BH_B$ and $CH_C$ intersect at one point.}

\includegraphics[scale=1]{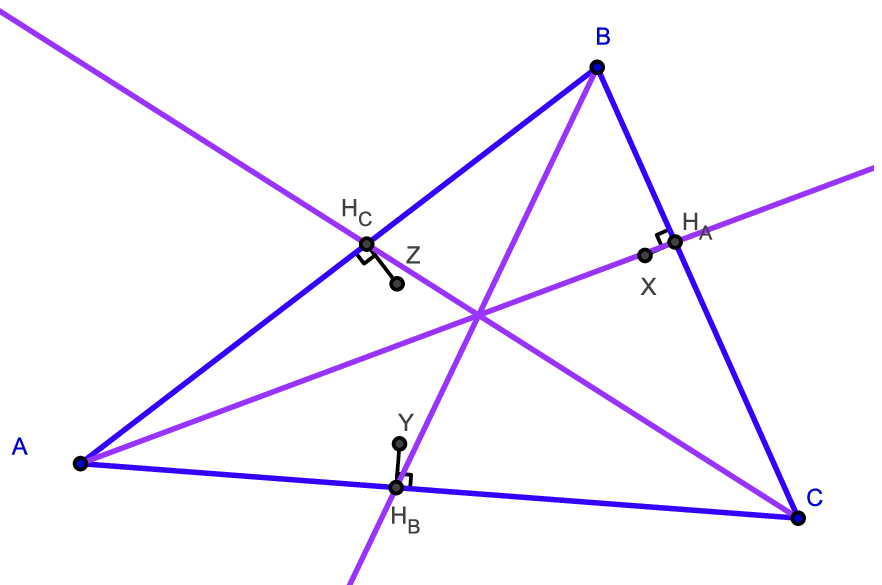}

\smallskip
\textbf {Remark.} Points $H'_A$, $H'_B$, $H'_C$ are defined similarly to points $H_A$, $H_B$ and $H_C$. 
Then points $H_A$, $H_B$,$H_C$, $H'_A$, $H'_B$ and $H'_C$ lie on the same conic section. 
This follows by a criterion for being conconic (see \S2), theorem 2 and Ceva's theorem.

\smallskip
\textbf {Theorem 3.} {\it Let $R'$ be the intersection point of $AX'$, $BY'$ and $CZ'$ 
(these lines intersect at a point by Theorem 1). 
Lines $XX'$, $YY'$, $ZZ'$ and $RR'$ are intersect at a point.}

\smallskip
\textbf {Remark.} Clearly, $X'$ is the isogonal conjugate of $X$, $Y'$ is the isogonal conjugate of $Y$ and 
$Z'$ is the isogonal conjugate of $Z'$.
This means that $R'$ is the isogonal conjugate of $R$.

\smallskip
Denote by $Q$ the intersection point of $XX'$, $YY'$, $ZZ'$ and $RR'$. 

\includegraphics[scale=1]{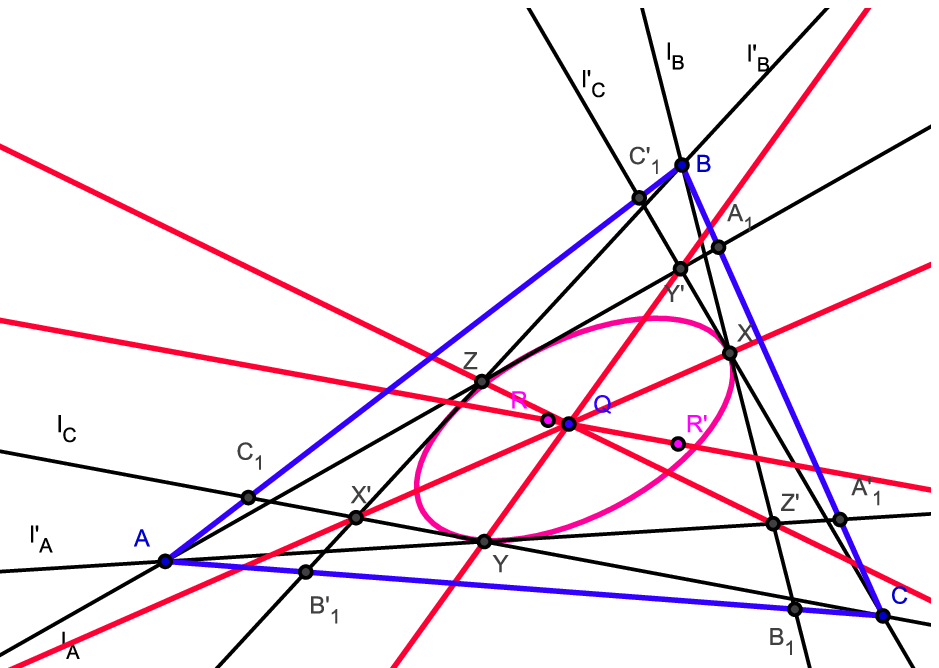}

\smallskip
\textbf{Theorem 4.} {\it Points $A_1$, $A'_1$, $B_1$, $B'_1$, $C_1$ and $C'_1$ lie on the same conic.}

\includegraphics{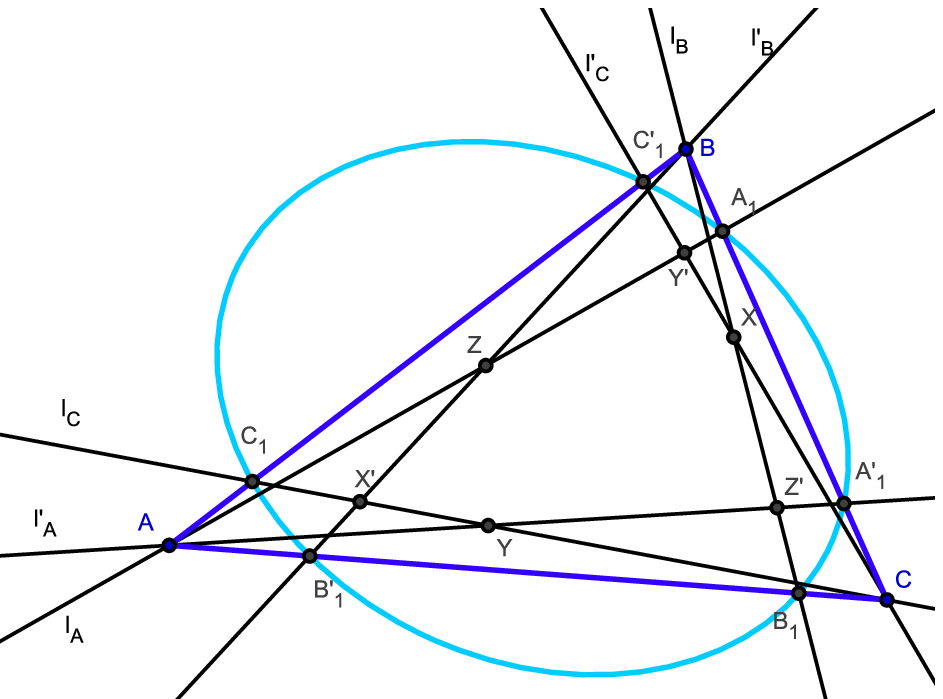}

\smallskip
A simple version of theorems is presented at the appendix.

\begin{center}
Proofs and some new corollaries
\end{center}

  {\it Proof of Theorem 1.} Cycles $\omega (\triangle XBC)$, $\omega (\triangle YAC)$ and $\omega (\triangle ZAB)$ are drawn through $\triangle XBC$, $\triangle YAC$ and $\triangle ZAB.$ Point $K$ is second point of intersection $AX$ and $\omega (\triangle XBC).$ Points $L$ and $M$ are similarly defined (for $BY$ and $CZ$ accordingly).

\includegraphics[scale=0.2]{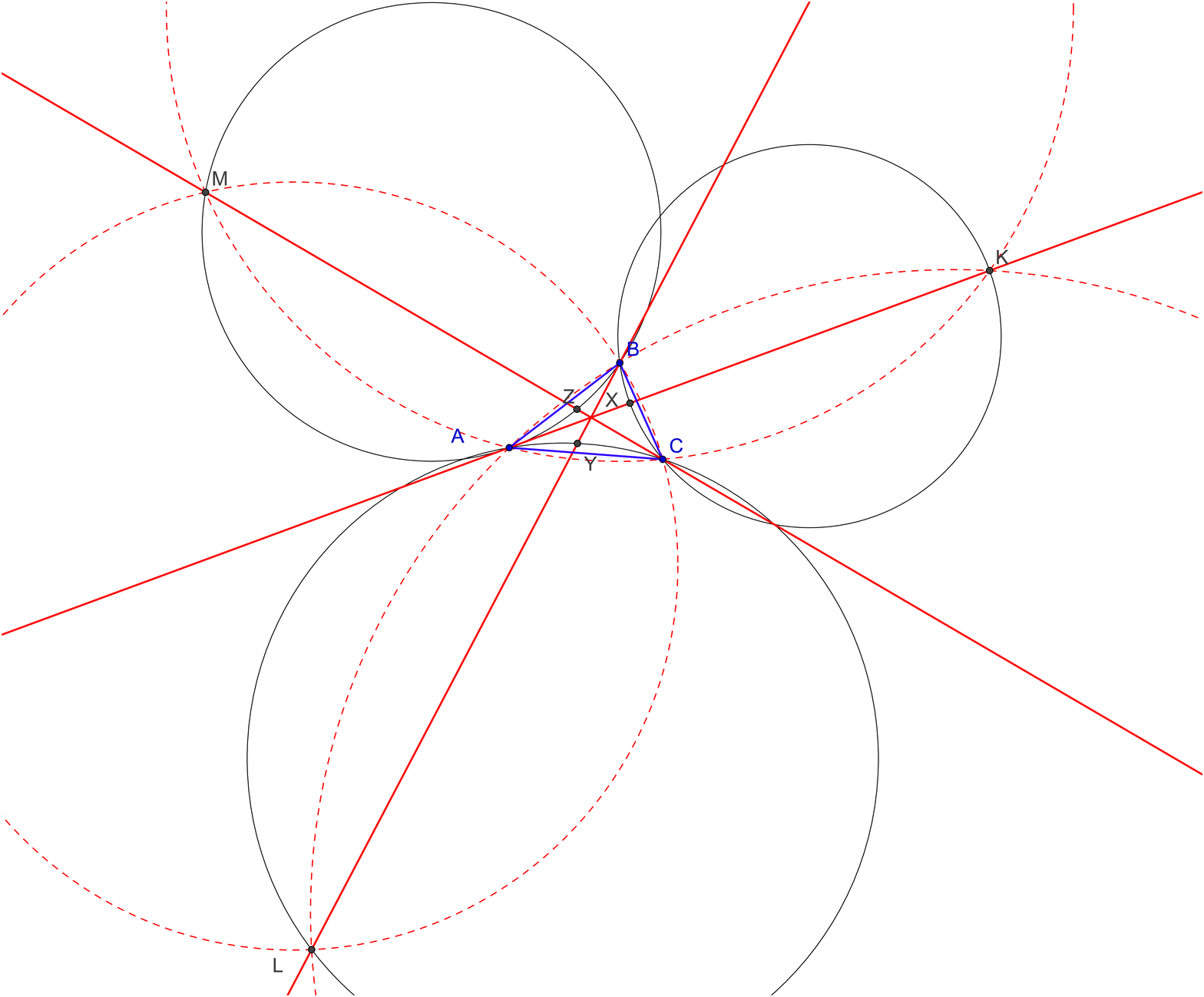}

  $$\angle CMB = \angle ZAB = \angle YAC = \angle BLC$$

We have that points $M,$ $L$, $B$ and $C$ are concyclic\footnote{This formula is true only for this picture (it means that arrangement of lines $l_A$, $l_B$ and $l_C$ can be very different), but in that cases we'll proof that points $M,$ $L$, $B$ and $C$ are concyclic in a similar way}. Let's define this cycle $\omega _x.$ Cycles $\omega _y$ and $\omega _z$ are similarly defined. It's clear that lines $AX$, $BY$ and $CZ$ are radical axises this cycles. This mean that $R$ (point of intersection lines $AX$, $BY$ and $CZ$) is radical centre.$\vartriangleleft$

 \smallskip
 {\it Proof of Corollary.} Let's fix lines $l_A$, $l_B$, $l_C$, $l'_A$, $l'_B.$ Consider conic, which this five lines tangent it (certainly, we can find the only one conic, which five lines general provision tangent this conic - [1]). This conic is inscribed in to $\angle ZAZ'$ and in to $\angle XBX'$. It means that point $P'$ (points $P$ and $P'$ are two focuses of this conic) is isogonal conjugate of $P$\footnote{Here I used isogonal properties of conics - [1]}. Consequently, tangent lines, which are drawn through vertice $C$, are isogonal lines.

{\it Proposition}: Let us fix $\triangle ABC$ and a conic. Tangent lines $d_{A}$, $d'_{A}$, $d_{B}$, $d'_{B}$, $d_{C}$ and $d'_{C}$ (this lines were drawn through vertices $A$, $B$ and $C$, respectively)  to conic create hexagon $XY'ZX'YZ'$. Lines $AX$, $BY$, $CZ$ are concur. This  proposition is conclusion of Brianchon's theorem\footnote{I used here next formulation of this theorem: Lines $l_i$, $i=1,...,6$ tangent the same conic section. Point $A_{ij}$ is point of intersection lines $l_i$, $l_j$. Lines $A_{12}A_{45}$, $A_{23}A_{56}$ and $A_{34}A_{61}$ are concur then.}  in next order: $d_A$, $d'_A$, $d_B$, $d'_B$, $d_C$, $d'_C$.

\smallskip
{\it Proof of Theorem 2.}

$$\frac{XH_A}{BH_A} = tg( \angle (BX,BC)) \Rightarrow \frac{BH_A}{CH_A} = \frac{tg( \angle (XC,BC))}{tg( \angle (BX,BC))}.$$

We can take that:
$$ \frac{CH_B}{AH_B} = \frac{tg( \angle (YA,AC))}{tg( \angle (YC,AC) )} = \frac{tg( \angle (AY,AC))}{tg( \angle (XC,BC))} ; $$
$$\frac{AH_C}{BH_C} = \frac{tg( \angle (BZ,AB))}{tg( \angle (AZ,AB) )} = \frac{tg( \angle (BX,BC))}{tg( \angle (AY,AC))} \Rightarrow$$%

$$\frac{AH_C}{BH_C} \cdot \frac{BH_A}{CH_A} \cdot \frac{CH_B}{AH_B} = 1 \Rightarrow $$
This formula shows that three lines are concur (there is converse proposition of Ceva's theorem).$\vartriangleleft$

\smallskip
{\it Proof of Theorem 1.} Let's proof general result: Let us fix $\triangle ABC$ and a conic. Tangent lines $d_{A}$, $d'_{A}$, $d_{B}$, $d'_{B}$, $d_{C}$ and $d'_{C}$ (this lines were drawn through vertices $A$, $B$ and $C$, respectively)  to conic create hexagon $XY'ZX'YZ'$. Points $R$ and $R'$ are defined like before. Then lines $XX'$, $YY'$, $ZZ'$ and $RR'$ are concur.

Let's proof next lemma:

\smallskip
\textbf {Lemma.} {\it Lines $XY$,$X'Y'$ and $AB$ are concur. Let's define this point $C'$. Points $A'$ and $B'$ are similarly defined.}

\includegraphics[scale=0.6]{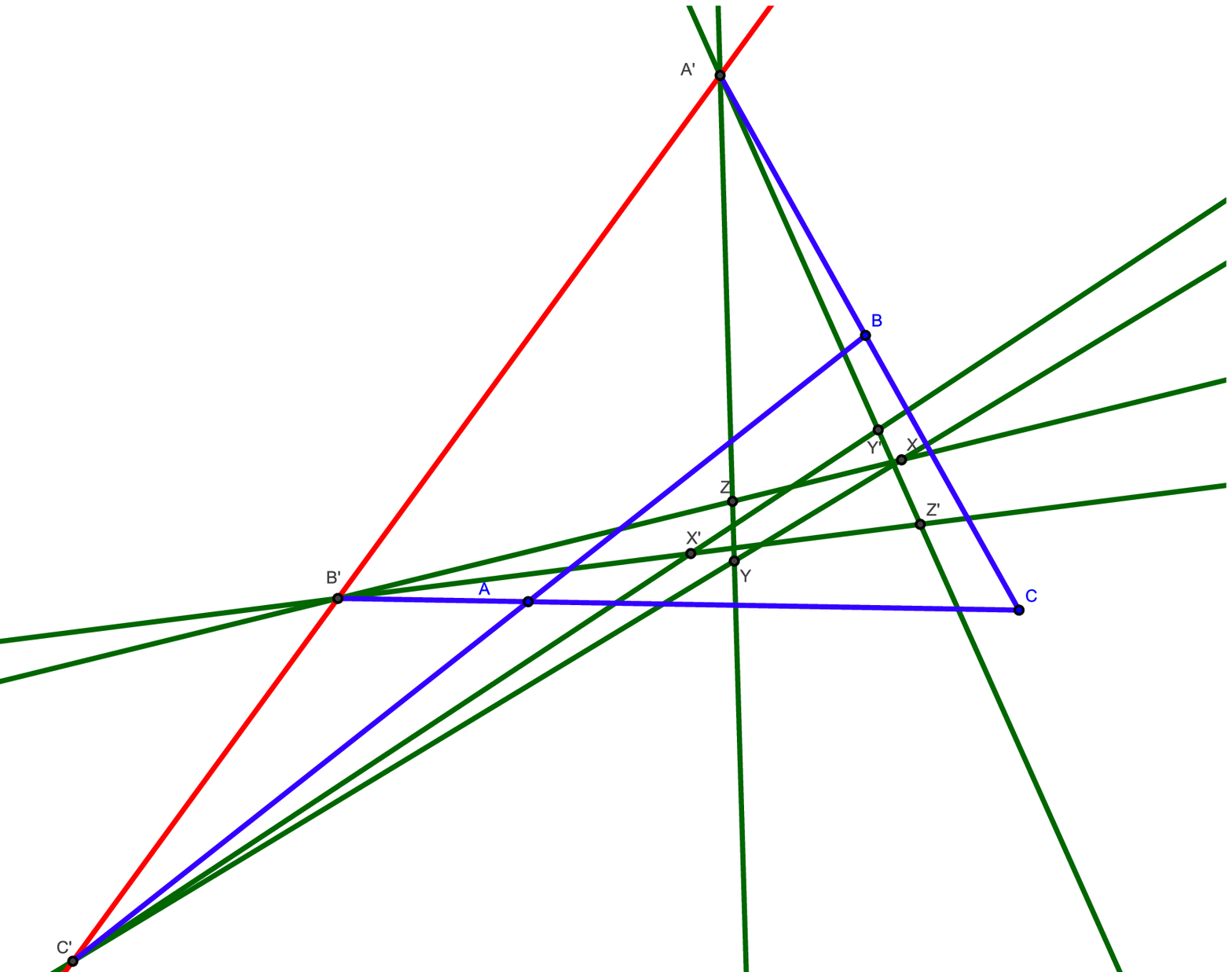}

{\it Proof.}: Let's look at $\triangle X'YB$ and $\triangle XY'A$: $X'Y$ cross with $XY'$ at a vertice $C$; $BY$ cross with $AX$ at a point $R$ and $X'B$ cross with $Y'A$ at a point $Z$. We have that points $C$, $R$ and $Z$ are collinear. It means  lines $XY$,$X'Y'$ and $AB$ are concur (there is converse proposition of Desargues'\footnote{If the three straight lines joining the corresponding vertices of two triangles $\triangle ABC$ and $\triangle A'B'C'$ all meet in a point (the perspector), then the three intersections of pairs of corresponding sides lie on a straight line (the perspectrix).} theorem).

Let's return to the proof of this theorem: we have that $\triangle XYZ$ and $\triangle ABC$ are perspective. This means that points $A'$, $B'$ and $C'$ are collinear (the perspectrix). Now we can conclude that $\triangle XYZ$ and $\triangle X'Y'Z'$ are perspective. It means that lines $XX'$, $YY'$ and $ZZ'$ are concur (but we can proof this propositision with Brianchon's theorem). $Q$ is the point of intersection $XX'$, $YY'$ and $ZZ'$. Let's look at $\triangle RYZ$ and $\triangle R'Y'Z'$: $RY$ cross with $R'Y'$ at  a vertice $B$, $RZ$ cross with $R'Z'$ at  a vertice $C$ and  $YZ$ cross with $Y'Z'$ at a point $A'$ (in lemma we proofed that $A'$ lays on $BC$). It means that $\triangle RYZ$ and $\triangle R'Y'Z'$ are perspective (there is converse proposition of Desargues' theorem). Now we have that point $Q$ lays on line $RR'$.$\vartriangleleft$

\smallskip
{\it Proof of Theorem 4.} This theorem is equivalent to the converse preposition which we generalized: Points $A_1$, $A'_1$, $B_1$, $B'_1$, $C_1$ and $C'_1$ lie on sidelines $BC$, $AC$ and $AB$ accordingly. Lines $d_{A}$, $d'_{A}$, $d_{B}$, $d'_{B}$, $d_{C}$ and $l'_{C}$ tangent the one conic if and only if $A_1$, $A'_1$, $B_1$, $B'_1$, $C_1$ and $C'_1$ lie on the same conic.

I used next proposition (task 14 in [1]) - criterion of conconic:

Let's fix $\triangle ABC.$ Points $A_1$, $A_2$ lie on sideline $BC$, $B_1$, $B_2$ lie on sideline $AC$ and $C_1$, $C_2$ lie on sideline $AB$. This six points are conconic if and only if, when
$$\frac{BA_1 \cdot BA_2}{CA_1 \cdot CA_2} \cdot \frac{CB_1 \cdot CB_2}{AB_1 \cdot AB_2} \cdot \frac{AC_1 \cdot AC_2}{BC_1 \cdot BC_2} = 1$$

Let $K$, $M$, $L$ be intersection points of lines $AX$ with $BC$, $BY$ with $AC$, $CZ$ with $AB$ accordingly. We can take (Ceva's equality):
$$\frac{CK}{BK} \frac{BC_2}{AC_2} \frac{AB_1}{CB_1} = 1; \frac{AM}{CM} \frac{CA_2}{BA_2} \frac{BC_1}{AC_1} = 1; \frac{BL}{AL} \frac{AB_2}{CB_2} \frac{CA_1}{BA_1} = 1$$

Let's multiply this three equalities and let's use criterion of conconic:
$$\frac{CK}{BK} \frac{AM}{CM} \frac{BL}{AL} = 1$$

We have that lines $AX$, $BY$ and $CZ$ are concur (there is converse proposition of Ceva's theorem).
Now, the proposition of this theorem is obvious (there is converse proposition of Brianchon's theorem or we can thinking like in theorem 1).$\vartriangleleft$

\smallskip
\textbf{Corollary 1.} Let $A_2$ be the intersection point of lines $C_1B_1$ and $C'_1B'_1$. Points $B_2$ and $C_2$ are similarly defined. We can have that point $A_2$ lie on line $XX'$ (Pappus's\footnote{Points $A_1$, $B_1$, $C_1$ are collinear and points $A_2$, $B_2$, $C_2$ are collinear too. Then three intersection points of lines $A_1B_2$ and $A_2B_1$, $B_1C_2$ and $B_2C_1$, $C_1A_2$ and $C_2A_1$ are incident to a (third) straight line (next Pappus's line).} theorem). We can receive next result after it:

Points $A_1$, $A_2$, $B_1$, $B_2$, $C_1$ and $C_2$ lie on sidelines $BC$, $AC$ and $AB$ accordingly. Three Pappus's lines are concur if and only if  points $A_1$, $A_2$, $B_1$, $B_2$, $C_1$ and $C_2$ are conconic.

\includegraphics{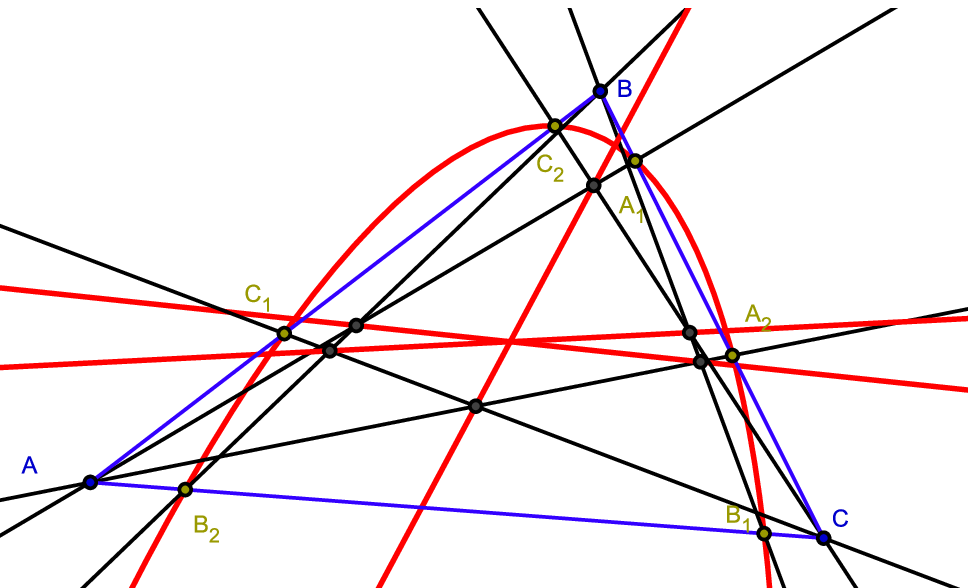}

\smallskip
\textbf{Corollary 2.} Let $A_3$ be the intersection point of $A_1C'_1$ and $B_1A'_1$. Points $B_3$ and $C_3$ are similarly defined. Then lines $AA_3$, $BB_3$ and $CC_3$ are concur.

{\it Proof.}: Pascal's\footnote{The theorem states that if a hexagon is inscribed in a conic, then the three points at which the pairs of opposite sides meet, lie on a straight line.

It is obvious that this theorem is generalization of Pappus's theorem.} theorem says that three points of intersection lines $AB$ with $A'_1B_1$, $BC$ with $B'_1C_1$ and $AC$ with $A_1C'_1$ are collinear. It means that $\triangle ABC$ and $\triangle A_3B_3C_3$ are perspective (there is converse proposition of Desargues' theorem).$\vartriangleleft$

\textbf{Remark.} We can this result generalize: Let's look at two triangles: $\triangle ABC$ and $\triangle A'B'C'.$ Points $C_1$ and $B'_1$ are points of intersection lines $C'B'$ with $AB$ and $AC$ accordingly. Points $A_1$, $A'_1$, $B_1$, $C'_1$ are similarly defined. This two triangles are perspective if and only if points $A_1$, $A'_1$, $B_1$, $B'_1$, $C_1$ and $C'_1$ are conconic\footnote{In case, when vertices first triangle lie on corresponding sidelines of second triangle, this proposition isn't true.}.

\includegraphics[scale=0.9]{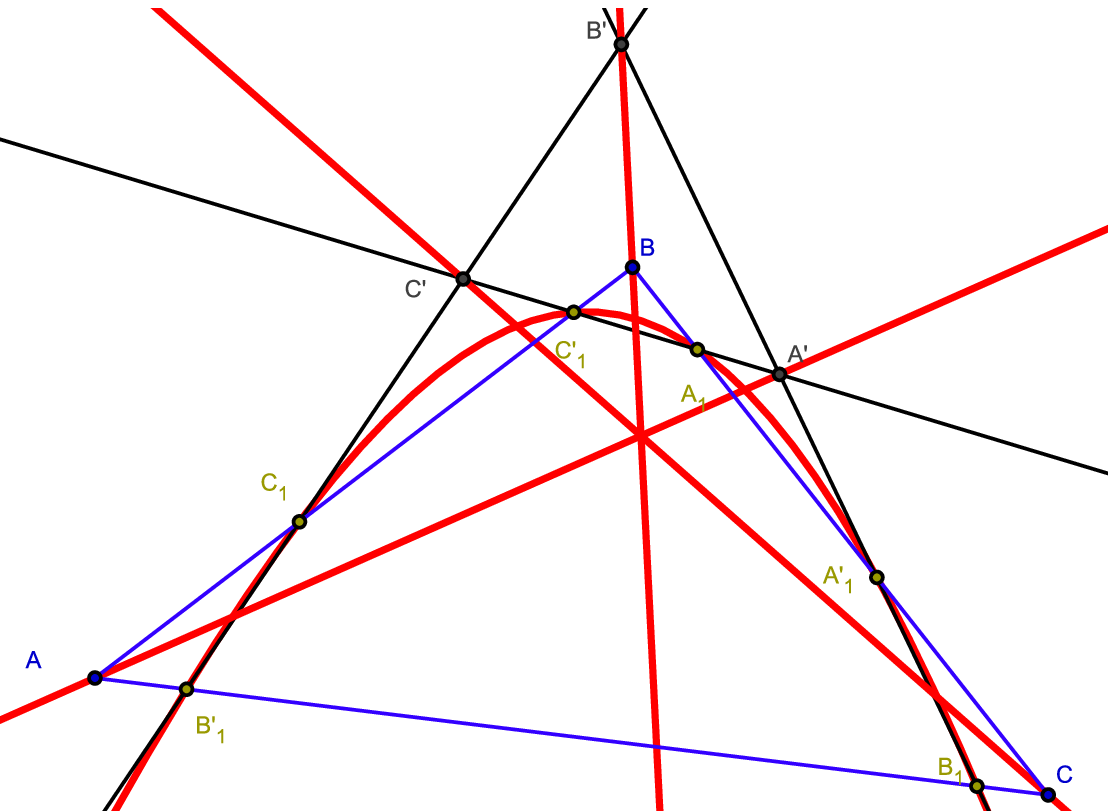}

\begin{center}
  A generalization of the Morley Point
\end{center}

\textbf{Morley's theorem}:  The three points of intersection of the adjacent trisectors of the angles of any triangle form an equilateral triangle.

And one theorem says that this two triangles are perspective. And point of intersection was named ${\it 2^{nd} Morley Centre} - X(357)$\footnote{This notation is notation of Clark Kimberling's Encyclopedia of Triangle Centers (ETC)}. We'll try to generalize this theorem next.

\includegraphics{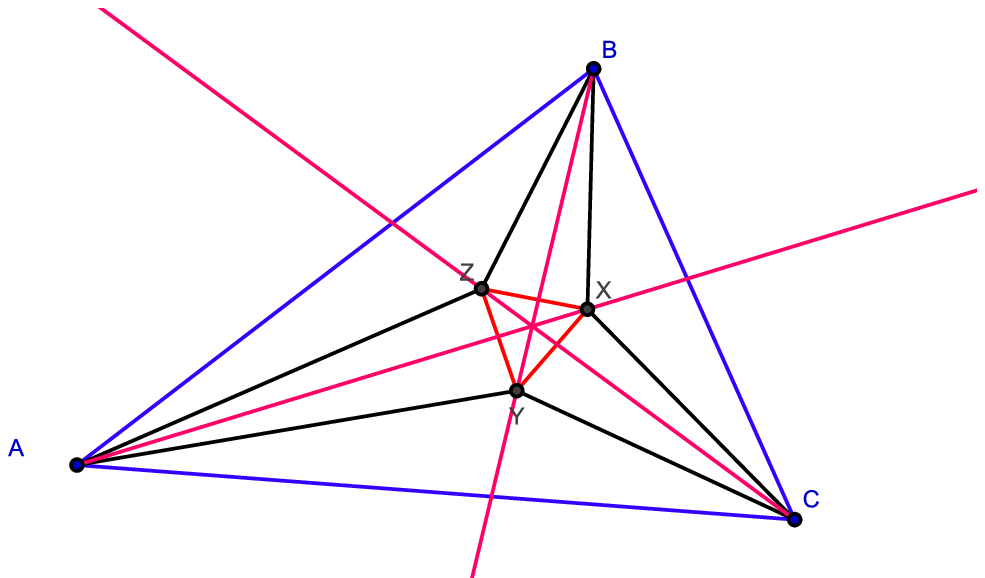}

Let's fix $k \in [-1;1]$. We will use the same notation in this construction like before, but with next difference:
$$\frac{\angle A_1AB}{\angle A} = \frac{\angle B_1BC}{\angle B} = \frac{\angle C_1CA}{\angle C} = k.$$

$k$ is positive in case, when lines $l_{A}$, $l_{B}$ and $l_{C}$ are drawn inside and $k$ is negative, when lines $l_{A}$, $l_{B}$ and $l_{C}$ are drawn outside.

It is clear that lines $AX$, $BY$ and $CZ$ intersect at a point (for $k=1/3$ the point of intersection is ${\it 2^{nd} Morley Centre}$). Let's define the point of intersection $R(k).$ It is easy to check (maybe with Ceva's theorem) that point $R(k)$ has next barycentric coordinates:

$$R(k) = \left( \frac{sin (\angle A) \cdot sin(\angle A \cdot k)}{sin(\angle A \cdot (1-k))}: \frac{sin (\angle B) \cdot sin(\angle B \cdot k)}{sin(\angle B \cdot (1-k))}: \frac{sin (\angle C) \cdot sin(\angle C \cdot k)}{sin(\angle C \cdot (1-k))} \right)$$

It is easy to see (in this formula or in this geometrical construction) that the curve of the concurrence is self-isogonal curve. It means that $X(358)$ (isogonal conjugate of ${\it2^{nd} Morley Centre}$) lays on this curve. Incentre lays on this curve too.

When $k=-1:$ $R(-1) = -\frac{1}{2} (tg \angle A: tg \angle B: tg \angle C) = H.$ It means that $Circumcenter$ lays on this curve too, because $Circumcenter$ is isogonal conjugate of $Orthocenter.$

Let's find limit points:

$k=0:$ In this case $$\frac{AL}{BL} = lim_{k \mapsto 0} \frac{sin(\angle B) \cdot sin(\angle A \cdot (1-k)) \cdot sin (\angle B \cdot k)}{sin(\angle A) \cdot sin(\angle B \cdot (1-k)) \cdot sin (\angle A \cdot k)} = \frac{sin(\angle B \cdot k)}{sin(\angle A \cdot k)} = \frac{\angle B}{\angle A}$$

Here I used that $lim_{x \mapsto 0} \frac{sinx}{x}=1.$ We have

$$R(0) = (\angle A: \angle B: \angle C)$$

$k=1:$ $$\frac{AL}{BL} = lim_{k \mapsto 1} \frac{sin(\angle B) \cdot sin(\angle A \cdot (1-k)) \cdot sin (\angle B \cdot k)}{sin(\angle A) \cdot sin(\angle B \cdot (1-k)) \cdot sin (\angle A \cdot k)} =$$

$$= \frac{sin^2 (\angle B) \cdot sin(\angle A \cdot (1-k))}{sin^2 (\angle A) \cdot sin(\angle B \cdot (1-k))} = \frac{sin^2 (\angle B) \cdot \angle A}{sin^2 (\angle A) \cdot \angle B}$$

Here I used that $lim_{x \mapsto 0} \frac{sinx}{x}=1$ too (like before). We have that

$$R(1) = \left (\frac{sin^2 \angle A}{\angle A}: \frac{sin^2 \angle B}{\angle B}: \frac{sin^2 \angle C}{\angle C} \right )$$

It means that $R(0)$ is isogonal conjugate of $R(1).$

It is clear that all results in general construction (\textbf{theorem 1-4}) are true here too. Let's look at the special case of \textbf{theorem 2}. Let's define this point of intersection $D(k).$ In case when $k=1/2$ we can receive $GergonnePoint.$ Easy to find that $$D(k)=(tg(\angle A \cdot k): tg(\angle B \cdot k): tg(\angle C \cdot k)).$$ When $k=1$, $D(1)=H.$ In a similar way we can take that $D(0) = R(0)$, like in the foregoing.

There is very difficult to find barycentric coordinates of $Q(k)$, but it is not difficult to say her limit points: $Q(0) = R(0)$ and $Q(1)=R(1).$

\begin{center}
  Obverse
\end{center}

Lines $l_A$ and $l'_A$ are isogonal lines in general design. In this construction we'll consider case where lines $l_A$ and $l'_A$ are isotomic lines ($l_B$ and $l'_B$, $l_C$ and $l'_C$ are isotomic lines too next).

Some notation: $\triangle A_0B_0C_0$ is the median triangle.

\includegraphics{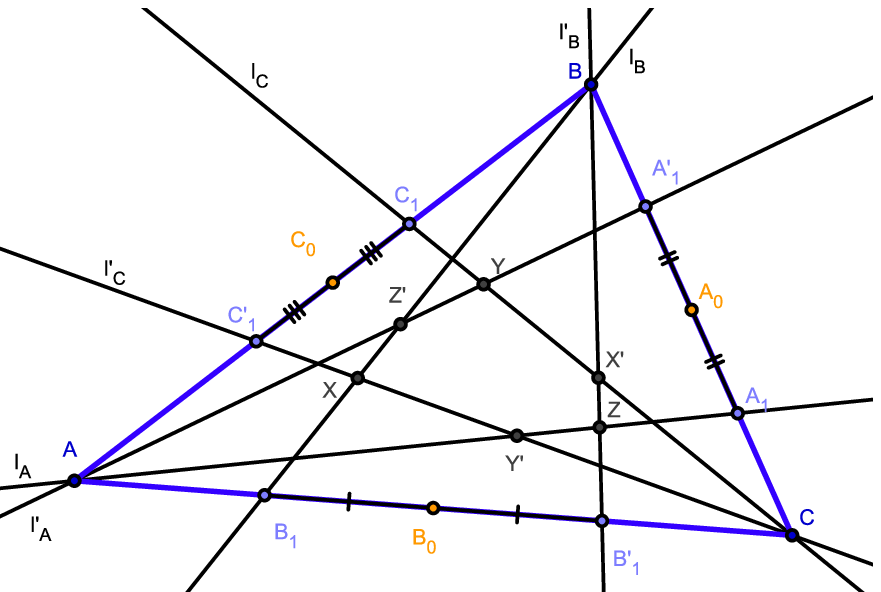}

Then this propositions are true:

\smallskip
1. Lines $l_{A}$, $l'_{A}$, $l_{B}$, $l'_{B}$, $l_{C}$ and $l'_{C}$ tangent the one conic. We can conclude that lines $AX$, $BY$ and $CZ$ intersect at a point ($P$) then. Properly lines $AX'$, $BY'$ and $CZ'$ are concur too ($P'$ is the point of intersection). And $P'$ is isotomic conjugate of $P.$

\smallskip
2. Lines $XX'$, $YY'$, $ZZ'$ and $PP'$ are concur.

\smallskip	
3. Points $A_1$, $A'_1$, $B_1$, $B'_1$, $C_1$ and $C'_1$ lie on the same conic section.

\smallskip
4.	Point $A_2$ is point of intersection $A_1C'_1$ and $B_1A'_1.$ Points $B_2$ and $C_2$ are similarly defined. Lines $AA_2$, $BB_2$ and $CC_2$ are concur then.

{\it Proof.}: Let's make affine transformation in which $\triangle ABC$ will translate in equilateral triangle. Thus points $A_0$, $B_0$ and $C_0$ will pass in the middle of a new triangle (because affine transformations save attitude of parallel pieces). Points $A_1$, $B_1$ and $C_1$ will pass in new points, but they will remain to lie on the corresponding edges of a new triangle. And points $A'_1$, $B'_1$ and $C'_1$ will pass in the points which they are symmetric to $A_1$, $B_1$ and $C_1$ rather the middle of the corresponding sides (besides because affine transformation save attitudes on a line). It means that the condition of theorems didn't change after transformation.

As points $A_1$ and $A'_1$ are symmetric rather $A_0$ and as a $\triangle ABC$ is equilateral triangle, then $\angle A_1AA_0 = \angle S'_1AA_0.$ We can receive that $l'_A$ is symmetric $l_A$ with respect to bisecting line $AA_0$. Now we have construction which is analogously with general construction (but only for equilateral triangle).$\vartriangleleft$

Let's explore the special interesting case when pieces $A_1A'_1$, $B_1B'_1$ and $C_1C'_1$ are proportional to corresponding edges:
$$\frac{A_1A'_1}{BC}=\frac{B_1B'_1}{AC}=\frac{C_1C'_1}{AB}.$$
It is clear (maybe with theorem about 4 points of trapezium or with affine transformations) that points $P$ and $P'$ will be $Centroid$ of $\triangle ABC.$  Analogously, points of intersection in propositions \textbf{2} and \textbf{4} will be  $Centroid$ of $\triangle ABC$ too.

\includegraphics{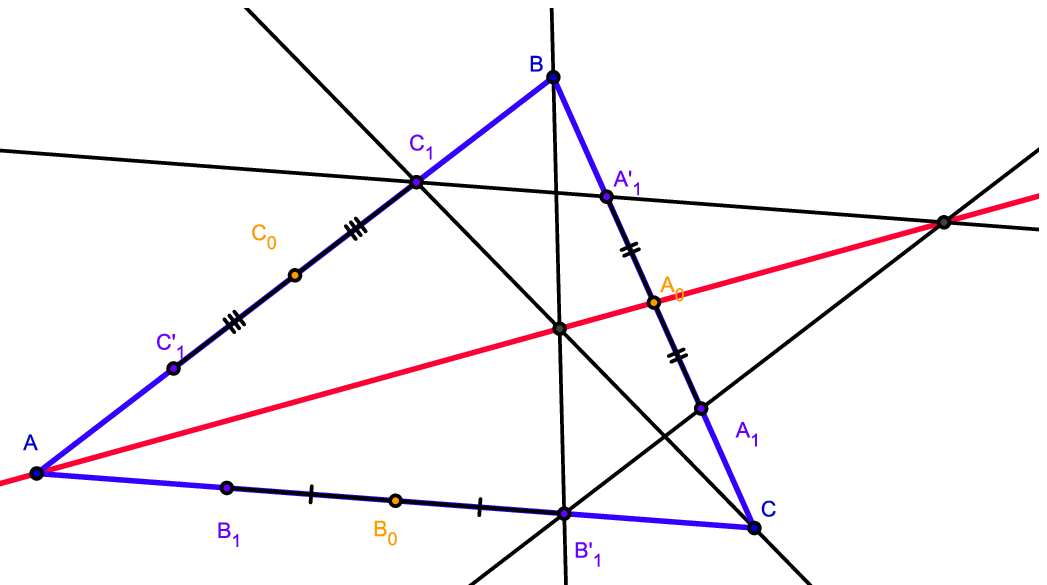}

\begin{center}
References
\end{center}

[1]A.Akopyan, А.Zaslavsky "Geometric properties of second order curves", М, MCCME, 2007.
\end{document}